\input amstex
\input amsppt.sty

\def\ni\noindent
\def\sbs{\subset}

\def\as{\operatorname{asdim}}
\def\diam{\operatorname{diam}}
\def\holim{\operatorname{holim}}
\def\dist{\operatorname{dist}}
\def\ker{\operatorname{ker}}
\def\cof{\operatorname{cof}}
\def\Map{\operatorname{Map}}
\def\Lip{\operatorname{Lip}}
\def\Tor{\operatorname{Tor}}
\def\Nerve{\operatorname{Nerve}}
\def\R{\text{\bf R}}

\def\Q{\text{\bf Q}}
\def\Z{\text{\bf Z}}
\def\H{\text{\bf H}}
\def\L{\text{\bf L}}

\def\E{\text{\bf E}}

\def\bS{\text{\bf S}}
\def\BG{\text{\rm B}\Gamma}
\def\EG{\text{\rm E}\Gamma}

\def\sA{\Cal A}

\def\sU{\Cal U}

\def\t{\tilde}

\def\x{\times}

\def\mk{\vskip .1in}

\hoffset= .8in
\voffset= 1.2in
\hsize=5.0in
\vsize=6.5in
\baselineskip=24pt
\NoBlackBoxes
\topmatter
\author
A.N. Dranishnikov, S. Ferry, and S, Weinberger
\endauthor

\title
An Etale Approach to the Novikov Conjecture
\endtitle
\abstract We show that the rational Novikov conjecture for a group
$\Gamma$ of finite homological type follows from the mod 2
acyclicity of the Higson compactifcation of an E$\Gamma$. We then
show that for groups of finite asymptotic dimension the Higson
compactification is mod p acyclic for all p, and deduce the integral
Novikov conjecture for these groups.
\endabstract

\thanks This work was partially supported by supported by the NSF.
\endthanks

\address University of Florida, Department of Mathematics, P.O.~Box~118105,
358 Little Hall, Gainesville, FL 32611-8105, USA
\endaddress

\address Department of Mathematics, Rutgers University, 110
Frelinghuysen Road,
Piscataway, NJ 08854-8019
\endaddress

\address Department of Mathematics, The
University of Chicago, Chicago, IL 60637
\endaddress


\email  dranish\@math.ufl.edu
\endemail

\endtopmatter

\document
\head \S1 Introduction \endhead

\mk\noindent Ten years ago, the most popular approach to the Novikov conjecture went
via compactifications. If a compact aspherical manifold, say, has a universal cover
which suitably equivariantly compactifies, already Farrell and Hsiang [FH] proved
that the Novikov conjecture follows. Subsequent work by many authors
weakened their hypotheses and extended the idea to other settings. (See
e.g. [CP1, CP2, FW, Ro1].)

\mk\noindent In recent years other coarse methods have supplanted
the compactification method (most notably the embedding method of
[STY] in the C* algebra setting). The reason for this is that it
seemed to have a better chance of applying generally, while
compactifications effective for proving the Novikov conjecture
seemed to require some special geometry for their construction. For
a brief time, it seemed that this could conceivably not be the case.
Higson introduced a general compactification of metric spaces
(somewhat reminiscent of the Stone-\v Cech compactification) that
automatically has half of the properties necessary for application
to the Novikov conjecture. The missing property was acyclicity,
which holds for the Stone \v Cech compactification in dimensions
$>1$, so it was natural to hope that the Higson compactification of
a universal cover is automatically acyclic. (See [Ro1].)

\mk\noindent Unfortunately, it was soon realized by 
and others (see [Kee, DF])
that the Higson compactification, even for manifolds as small as \R, has nontrivial
rational cohomology. Thus, it was felt that general compactifications were not
suitable for the problem -- one has to use geometric compactifications.
However, Gromov has recently harpooned the embedding approach by
constructing finitely generated groups which do not uniformly embed in Hilbert
space [G2]. Moreover, a number of authors (e.g. [HLS] [O]) have shown that
GromovÍs groups can be used to construct counterexamples to general forms of
the Baum-Connes conjecture.

\mk\noindent Moreover, for reasons that are not entirely clear, the embedding method
has never been translated into pure topology; the results on the integral Novikov
conjecture so obtained, from the L-theoretic viewpoint never account for the
prime 2, and there do not seem to be many results in pure algebraic K-theory (or
A-theory) provable by this method.

\mk\noindent This paper seeks to rehabilitate the Higson
compactification approach. We will show that

\proclaim{Theorem 1} If the Higson compactification of $E\Gamma$ is mod 2
acyclic, and $B\Gamma$ is finite type, then the  Novikov
conjecture for $\Gamma$ holds at the prime 2.  \endproclaim

\mk \noindent In other words, the assembly map localized at the prime 
two is an injection on homotopy groups. At odd primes we do not
know how to prove any corresponding statement, in general, for reasons
related to the examples in [DFW]: The L-spectrum away from 2 is
periodic K-theory, and cohomologically acyclic spaces may not be
acyclic for periodic K-theory.

\mk\noindent This difficulty is an infinite dimensional phenomenon, so for groups of
ñfinite asymptotic dimensionî this issue does not arise, and it is
possible to prove integral results.

\proclaim{Theorem 2} If $B\Gamma$ is a finite complex and $\Gamma$ has
finite asymptotic dimension and R\ is a ring with involution that has finitely 
generated L-groups, then the $\L(R)$ assembly map for $\Gamma$ is
injective. \endproclaim

\mk\noindent In this paper, by $\L(R)$, we will always mean $\L^{-\infty}(R)$.
\mk\noindent Theorem 2 is not a new theorem.  For $K(C^{*}\Gamma)$, this 
is due to [Yu1].
In algebraic K-theory there are papers [Ba, CG2], which apply to
L-theory as well (without the restrictions on R); moreover [CFY] has
yet another approach to the stronger L-theory result.  The novelty is
the method of using finite primes to repair the acyclicity of the
Higson compactification and then using this ``etale'' idea to get a
Novikov result.\footnote{Note that the celebrated paper of [BHM] on
the algebraic K-theory analogue of the Novikov conjecture is a very
nice illustration of a completely different etale idea. We thank Alain 
Connes for suggesting this evocative description of our method.}

\mk\noindent
Theorem 2 follows from the general discussion above, together with our third theorem.

\proclaim{Theorem 3} If $\Gamma$ is as in Theorem 2, then the
Higson compactification of $B\Gamma$ is mod p acyclic for every
$p$.\endproclaim

\mk\noindent To motivate and clarify our results, we close with a 
brief discussion 
of the Higson compactification and of the Keesling elements. If X is a metric
space, and $f: X \to \R$ is a bounded continuous function, we say that $f$ has
{\it decaying variation} if $$\lim \diam (f(B_{r}(x))) = 0$$
for any fixed r and any sequence of points x going to infinity in X.
Higson's
compactification is the smallest one containing $X$ densely so that all functions
with decaying variation extend. Let $C_{h}$ be the aggregate of such functions. We
obtain the Higson compactification of $X$ by embedding $X $ into 
$Maps[C_{h} : R]$ and taking the closure of the image.

\mk\noindent Let us now consider the first cohomology of the compactification of a
uniformly simply connected $X$. We represent $\check H^{1}(\bar X)$ by 
maps $f:\bar X \to S^{1}$. On $X$, 
we can lift $f|$ and view our class as being given by $\exp(2\pi i 
\tilde f)$, where
$\tilde f$ has decaying
variation but is not necessarily bounded. Maps $f,\,g:\bar X \to 
S^{1}$ represent the same class
iff $\tilde f$
and $\tilde g$ differ by a bounded function.

\mk\noindent Note therefore, that $H^{1}( \bar X; Z)$ is now an \R-vector space -- one can multiply
the function $\tilde f$ by a real number. In particular, this cohomology vanishes with mod
p coefficients. In this regard, the cohomology of the Higson compactification
resembles the uniformly finite homology of a uniformly contractible n
manifold, [BW], which is known to be an \R\  vector space, aside from the class
coming from the fundamental class of $X.$

\mk\noindent Our paper is organized as follows. Section 2 is devoted to ñdescentî, i.e.
deducing Novikov conjectures from metric results. The last section is devoted to
verifying the mod p acyclicity result for the finite asymptotic dimension 
case. This is
done using ideas of ñquantitative algebraic topologyî. The key idea is that when
homotopy groups are finite, it is sometimes possible to get extra Lipschitz
conditions on maps, a priori.

\head  \S2 $\L$-theory and assembly map with $\Z_p$
coefficients\endhead

\noindent For every ring  with involution $R$ (even for any
additive category with involution) Ranicki defines a 4-periodic spectrum
$\L_*(R)$ with $\pi_i(\L_*(R))=L_i(R)$ [Ran1,2]. We
use the notation $\L=\L_*(\Z)$. Strictly speaking, $\L$-spectra are 
indexed by $K$-groups. Our $\L(\R)$ equals $\L^{-\infty}(\R)$, the 
limit of the $\L$-spectra associated to the negative $K$-groups. This 
applies to all $\L$-spectra occurring in this paper.

\mk\noindent For a metric space $(X,d)$ we will denote by $C_X(R)$ the boundedly
controlled Pedersen-Weibel category whose objects are locally
finite direct sums $A=\oplus_{x\in X}A(x)$ of finite dimensional
free $R$-modules and whose morphisms are given by matrices with bounded
propagations (see [FRR]). For a subset $V\subset X$ we denote by $A(V)$ the
sum $\oplus_{x\in V}A(x)$. Ranicki [Ran2] defined {\it X-bounded
quadratic $L$-groups}\ \ $L_*(C_X(R))$ and the corresponding
spectrum $\L_*(C_X(R))$. We will use the notation
$\L^{bdd}(X)=\L_*(C_X(\Z))$.

\mk\noindent Suppose that $\bar X$ is a compactification of $X$ with corona
$Y=\bar X\setminus X$. Then one can define the {\it continuously controlled
category} $B_{X,Y}(R)$ by taking the same objects as above with
morphisms $f:A\to B$ in $B_{X,Y}(R)$ taken to be homomorphisms such that for every $y\in
Y$ and every neighborhood $y \in U\subset\bar X$ there is a smaller
neighborhood $y \in V\subset U$ such that $f(A(V))\subset A(U)$. This
category is additive and hence also admits an $L$-theory. The
corresponding spectrum for $R=\Z$ we denote by
$\L^{cc}(X)=\L_*(B_{X,Y}(\Z))$.

\mk\noindent If $\E=\{E_k\mid k\in\Z\}$ is a spectrum,  we will denote the
$\E$-homology of $X$ by $H_i(X;\E)$. If $X$ is a complex, then $H_i(X;\E)=\pi_i(X_+\wedge\E)
=\lim_{k\to\infty}\pi_{i+k}(X_+\wedge E_k)$. If $X$ is a compact metric space,
 we denote the Steenrod $\E$-homology of $X$ by $H_i(X;\E)$, i.e.
$H_i(X;\E)=\pi_i(\holim\{N^i_+\wedge \E\})$ where $X=\varprojlim N^i$
is the inverse limit of polyhedra [CP2], [EH].

\mk\noindent The following theorem was proven in [P], [CP1], [CP2].
 \proclaim{Theorem 2.1}
$L_i(B_{X,Y}(R))=\bar H_{i-1}(Y;\L(R))$ for all $i$ where the homology
on the right is reduced Steenrod $\L_*(R)$-homology.
\endproclaim

\mk\noindent A metric space $X$ is called {\it proper} if every ball $B_r(x)$
in $X$ is compact. Subsets $A,B\subset X$ of a proper metric space $(X,d)$ are called {\it
diverging} if $$\lim_{r\to\infty}d(A\setminus B_r(x_0),B\setminus
B_r(x_0))=\infty$$ where $x_0\in X$ is any basepoint, $B_r(x_0)$ is
the $r$-ball centered at $x_0$, and the distance $d(A',B')$
between sets is the infimum of distances between $d(a,b)$, $a\in
A'$, $b\in B'$. We note that the minimal compactification of a proper 
metric space $X$ with respect
to the property that the closures of every pair of diverging subsets in
$X$ have empty intersection in the corona is  the 
Higson compactification.

\mk\noindent A compactification $\tilde X$ of $X$ is called {\it Higson dominated}
if the identity map $1_X:X\to X$ admits a continuous extension
$\bar X\to\tilde X$, where $\bar X$ is the Higson compactification.
We note that for every Higson dominated
compactification $\t X$ of $X$ with corona $Y=\t X-X$, there is a forgetful
functor $C_X(R)\to B_{X,Y}(R)$. This functor defines a 
map of spectra $\phi:\L^{bdd}(X)\to\L^{cc}(X)$.

\mk\noindent If $X$ is the universal covering of a finite complex with
fundamental group $\Gamma$, then $\Gamma$ acts on $C_X(R)$ and
hence on $\L_*(C_X(R))$ with fixed set
 $\L_*(C_X(R))^{\Gamma}=\L_*(\text{\rm R}\Gamma)$ [CP].  If a compactification of $X$ with corona $Y$
is equivariant, then $\Gamma$ acts on $B_{X,Y}(R)$ and hence on
$\L_*(B_{X,Y}(R))$.

\mk\noindent We recall that the homotopy fixed set $X^{h\Gamma}$ of a pointed
space $X$ with a $\Gamma$ action on it is defined as the space of
equivariant maps $\Map_{\Gamma}(\EG_+,X)$.

\mk\noindent For a general spectrum $\E$ and a torsion free group $\Gamma$ it
was proven [C], [CP1] that
$$
\H_*(\BG;\E)\cong\H^{lf}_*(\EG;\E)^{\Gamma}
\cong\H^{lf}_*(\EG;\E)^{h\Gamma}.
$$
\mk\noindent  For a
locally compact space $X$, $H_i^{lf}(X;\E)$ is defined to be the
Steenrod $\E$-homology rel infinity of the one point compactification
$\alpha X$.  We note that for a complex $N$, $\H_*(N;\E)=N_+\wedge\E$
and for a compact space $Y$,
$\H_*(Y;\E)=\holim\{N^{\alpha}_+\wedge\E\}$ where $N^{\alpha}$ runs
over nerves of all finite open covers of $Y$.  If $\Gamma$ acts on a
compact space $Y$, then it acts on the set of all finite open covers
of $Y$ and hence on the spectrum $\H_*(Y;\E)$.

\mk\noindent The following theorem is due to 
Carlsson-Pedersen [CP1],[CP2]. See also Ranicki [Ran3]. The first part of it
is discussed in a more general setting in [Ros], Theorem 3.3.

\proclaim{Theorem 2.2}
For every group $\Gamma$ with finite classifying complex
$\BG$  there is a morphism of spectra called the
bounded control assembly map
$$A^{bdd}:\H^{lf}_*(\EG;\L)\to\L^{bdd}(\EG)$$ which fits into
a homotopy commutative diagram

$$
\CD
\H_*(\BG;\L) @ >A>>\L_*(\Z\Gamma)\\
@V{\simeq}VV @V{trf}VV\\
\H_*^{lf}(\EG;\L)^{h\Gamma} @ >A^{bdd,h\Gamma}>> \L^{bdd}(\EG)^{h\Gamma}\\
\endCD
$$
where $A$ is the standard assembly map, and the vertical arrows are
the natural maps from fixed sets to homotopy fixed sets.

\mk\noindent If the universal covering space $\EG$ admits a Higson
dominated equivariant compactification $X$ with corona $Y$
then the diagram can be extended
$$
\CD
\H_*(\BG;\L) @ >A>>\L_*(\Z\Gamma)\\
@V{\simeq}VV @V{trf}VV\\
\H_*^{lf}(\EG;\L)^{h\Gamma} @ >A^{bdd,h\Gamma}>> \L^{bdd}(\EG)^{h\Gamma}\\
@V=VV @V{\phi^{h\Gamma}}VV\\
\H_*^{lf}(\EG;\L)^{h\Gamma} @ >A^{cc,h\Gamma}>> \L^{cc}(\EG)^{h\Gamma}\\
\endCD
$$
where $A^{cc,h\Gamma}=\phi^{h\Gamma}\circ A^{bdd,h\Gamma}$. Moreover, $A^{cc}$ coincides
with the boundary map in the Steenrod $\L$-homology exact
sequence of pair $(X,Y)$:
$$A_p^{cc}=\partial :\H_*^{lf}(\EG;\L)\to \H_{*-1}(Y;\L).$$
\endproclaim

\mk\noindent Let $M(p)$ denote the Moore spectrum for the group 
$\Z_p=\Z/p\Z$.

\proclaim{Lemma 2.1}
Let $X$ be a finite-dimensional compact metric space.
Suppose that $X$ is $\Z_p$-acyclic, i.e.,
$\tilde H^*(X;\Z_p)=0$ for reduced \v{C}ech cohomology. Then
$X$ is acyclic for the reduced
Steenrod $\L\wedge M(p^k)$-homology
for all $k$.
\endproclaim
\demo{Proof}
In view of the universal coefficient formula [M], we note that $\tilde H_*(X;\Z_p)=0$ for reduced Steenrod 
homology.
The coefficient exact sequence shows that
$\tilde H_*(X;\Z_{p^k})=0$ for all $k$.
For every spectrum $\bS$ there is
a Steenrod homology Atiyah-Hirzebruch spectral sequence
with $$
E^2_{i,j}=\tilde H_i(X;H_j(S^0,\bS))
$$
which converges to $\tilde H_*(X;\bS)$ provided that $X$ is a finite dimensional
compact metric space [EH],[KS]. Since all of the groups $H_q(S^0,\L\wedge 
M(p^k))$ are $p$-primary, we have $E^2_{i,j}=0$ for all $i$ and $j$.\qed
\enddemo

\mk\noindent REMARK 1. The finite dimensionality condition is essential here.
There are acyclic compacta that have nontrivial
mod $p$ complex $K$-theory [T]. Namely, by results of Adams and Toda
there is a map $f:\Sigma^kM(\Z_p,m)\to M(\Z_p,m)$ which induces an isomorphism
in $K$-theory. The inverse limit $X$ of suspensions of $f$ is a
an acyclic compactum, since all bonding maps are trivial in cohomology,
yet it has nontrivial $K$-theory and nontrivial mod $p$ $K$-theory. 
Hence, for $p$ odd, $X$ has nontrivial mod $p$ $L$-theory.

\mk\noindent REMARK 2. The spectrum $\L$ localized at 2 is equivalent (see [MM]) to the
Eilenberg-MacLane 4-periodic spectrum generated by the space
$$
\prod_{i=1}^{\infty}K(\Z_{(2)},4i)\times K(\Z_2,4i-2).
$$
\mk\noindent Thus, $\L\wedge M(2)=\L_{(2)}\wedge M(2)$ is the Eilenberg-MacLane
spectrum generated by the space
$$
\prod_{i=1}^{\infty}K(\Z_{2},4i)\times K(\Z_2,4i-1)\times K(\Z_2,4i-2).
$$
\mk\noindent  Indeed, for any $R$, $\L(R)$ is Eilenberg-MacLane at 2 
[TW]. Thus, if a compactum $X$ is $\Z_2$-acyclic, then it is
$\L\wedge M(2)$-acyclic without any finite dimensionality assumption
on $X$.

\mk\noindent DEFINITION. The mod $p$ assembly map $A_p$ is 
$$A\wedge 1_{M(p)}:\H_*(\BG;\L)\wedge
M(p)\to\L_*(\Z\Gamma)\wedge M(p).
$$

\proclaim{Theorem 2.3} Suppose that the universal cover $X$ of a
finite aspherical complex $B$ has finite dimensional mod p acyclic Higson
compactification. Then the
mod $p^k$ assembly map for $\Gamma=\pi_1(B)$ is a split
monomorphism for every $k$.
\endproclaim
\demo{Proof} Let  let $m$ be the dimension of the Higson corona $\dim\nu X$. Using
Schepin's spectral theorem [Dr1] one can obtain a metrizable
$\Z_p$-acyclic $\Gamma$-equivariant compactification $\bar X$ of
$X$ with corona $Y$ such that $\dim Y=m$ (see the proof of Lemma 8.3.
in [Dr1]).
We introduce coefficients to the second
diagram of Theorem 2.2 by forming the smash product with $M(p^k)$.

$$
\CD
\H_*(B;\L)\wedge M(p^k) @ >A_p>>\L_*(\Z\Gamma)\wedge M(p^k)\\
@V{\simeq}VV @V{\phi^{h\Gamma}trf\wedge 1}VV\\
\H_*^{lf}(X;\L)^{h\Gamma}\wedge M(p^k) @ >A^{cc,h\Gamma}\wedge 1>>
\L^{cc}(X)^{h\Gamma}\wedge M(p^k)\\
\endCD
$$
Note also that $\L \wedge p^{k}=\cof(\x p^{k}:\L \to \L)$ and $\x 
p^{k}$ commutes with all of the structure mapsdefining the spectral 
structure and the $\Gamma$-action. 
We therefore have the commutative diagram
$$
\CD
\H_*(B;\L)\wedge M(p^k) @ >A_p>>\L_*(\Z\Gamma)\wedge M(p^k)\\
@V{\simeq}VV @V{\phi^{h\Gamma}trf\wedge 1}VV\\
(\H_*^{lf}(X;\L)\wedge M(p^k))^{h\Gamma} @>
(A^{cc}\wedge 1)^{h\Gamma}>> (\L^{cc}(X)\wedge M(p^k))^{h\Gamma}.\\
\endCD
$$
It suffices to show that 
$$
A^{cc}\wedge 1:\H_*^{lf}(X;\L)\wedge M(p^k)\to\L^{cc}(X)\wedge M(p^k)
$$
is an isomorphism because this implies that $(A^{cc}\wedge
1)^{h\Gamma}$ is an isomorphism and hence that $A_{p^{k}}$ is a split
monomorphism.

\mk\noindent In view of Theorems 2.1 and 2.2 we need to show that
$$
\partial\wedge 1:\H_*^{lf}(X;\L)\wedge M(p^k)\to
\H_{*-1}(Y;\L)\wedge M(p^k)
$$
is isomorphism. Note that $\partial\wedge 1_{M(p^k)}$ is
equivalent to the boundary homomorphism for the pair $(\bar X,Y)$
in $\L\wedge M(p^k)$-homology. By Lemma 2.1, $\H_*(\bar X;\L\wedge
M(p^k))=\H_*(pt;\L\wedge M(p^k))$ and hence $\partial\wedge 1_{M(p^k)}$ is
indeed an
isomorphism.
\qed
\enddemo
\mk\noindent It is known that the universal coefficient formula with $\Z_p$
coefficients (UCF) holds true for every generalized homology
theory and is natural with respect to morphisms of
spectra:

\proclaim{Proposition 2.1}
For every morphism of spectra $A:\E_1\to \E_2$ and every $p$ and $i$ there
is a commutative diagram
$$
\CD
0 @>>> \pi_i(\E_1)\otimes\Z_p @>>> \pi_i(\E_1\wedge M(p)) @>>>\pi_{i-1}(\E_1)\ast\Z_p @>>> 0\\
@. @V{A_*\otimes 1}VV @V(A\wedge 1)_*VV  @V(A_{*-1})\ast 1VV @.\\
0 @>>> \pi_i(\E_2)\otimes\Z_p @>>> \pi_i(\E_2\wedge M(p)) @>>>
\pi_{i-1}(\E_2)\ast\Z_p @>>> 0.\\
\endCD
$$
\endproclaim
\demo{Proof}
We apply the smash product with $\E_i$, $i=1,2$ to the cofibration of spectra
$\bS\to\bS\to M(p)$, where $\bS$ is the sphere spectrum.
Then the result follows from the homotopy exact sequence of the resulting
cofibrations
of spectra and the induced morphism between them
$$
\CD
E_1 @>>> E_1 @>>> E_1\wedge M(p)\\
@VVV @VVV @VVV\\
E_2 @>>> E_2 @>>> E_2\wedge M(p).\\
\endCD
$$
\qed
\enddemo

\proclaim{Theorem 2.4} Suppose that the  the universal cover of a finite 
complex $\BG$  has finite dimensional mod p acyclic Higson 
compactification. Then the integral assembly map $A$ is a monomorphism.
\endproclaim
\demo{Proof} In view of compactness of $\BG$ we have $\nu\Gamma=\nu \EG$.
By Theorem 2.3 $A\otimes 1_G$ is a split
monomorphism for every finite abelian group $G$. Since $\BG$
is a finite complex, the standard induction argument on the number
of cells show that the group $H_i(\BG;\L)$ is finitely generated for
every $i$. Hence for every $\alpha\in H_i(\BG;\L)$ there is $p$ such
that $\alpha\otimes 1\in H_i(\BG;\L)\otimes\Z_{p^{k}}$ is not zero.  By
the UCF there is a monomorphism $H_i(\BG;\L)\otimes\Z_{p^{k}}\to
H_i(\BG;\L(\Z)\wedge M(p^{k}))$ which, together with the assembly map, produces a
commutative diagram
$$
\CD
H_i(\BG;\L)\otimes\Z_{p^{k}} @>>> H_i(\BG;\L\wedge M(p^{k}))\\
@VA\otimes 1VV @VA_pVV\\
L_i(\Z\Gamma)\otimes\Z_{p^{k}} @>>> \pi_i(\L(\Z\Gamma)\wedge M(p^{k})).\\
\endCD
$$
This diagram implies that $A(\alpha)\ne 0$. \qed
\enddemo

\proclaim{Theorem 2.5} Suppose that $\Gamma$ is a group with $\BG$ a
finite complex.  If $\EG$ has mod 2 acyclic Higson compactification,
then the rational Novikov conjecture holds for
$\Gamma$.
\endproclaim
\demo{Proof}
In view of Remark 2 the argument of Theorem 2.3 for $p=2$ works without
the assumption $\dim \EG<\infty$. Hence the mod $2^k$ assembly map
$$
A\wedge 1_{M(2^k)}:\H_*(\BG;\L)\wedge
M(2^k)\to\L_*(\Z\Gamma)\wedge M(2^k)
$$
is a split monomorphism.
On the group level by the Universal Coefficient
Formula we have the diagram
$$
\CD
0 @>>> H_i(\BG;\L)\otimes\Z_{2^k} @>>> H_i(\BG;\L\wedge M(2^k)) @>>>
H_{i-1}(\BG;\L)\ast\Z_{2^k} @>>> 0\\
@. @V{A_*\otimes 1}VV @V(A\wedge 1)_*VV  @V(A_{*-1})\ast 1VV @.\\
0 @>>> \pi_i(\L_*(\Z\Gamma))\otimes\Z_{2^k} @>>>
\pi_i(\L_*(\Z\Gamma)\wedge M(2^k)) @>>>
\pi_{i-1}(\L_*(\Z\Gamma)\ast\Z_{2^k} @>>> 0.\\
\endCD
$$
Since the group $H_i(\BG;\L)$ is finitely
generated, it can be presented as $\oplus_{F_i}\Z\oplus \Tor_i$.
Since $A_*\otimes 1_{\Z_{2^k}}$ is a monomorphism, the kernel
$\ker(A_*)$ consists of $2^k$ divisible elements. Since $k$ is arbitrary,
$\ker(A_*|_{\oplus\Z})=0$. Therefore, $A_*\otimes 1_{\Q}$ is a monomorphism.
\qed
\enddemo

\noindent This proves Theorem 1 of the introduction.

\head \S3 Mod $p$ acyclicity of Higson compactifications of asymptotically finite
dimensional spaces\endhead

\mk\noindent We recall that a map $f:X\to Y$ between metric spaces is $\lambda$-Lipschitz
if $d_Y(f(x),f(x'))\le\lambda d_X(x,x')$ for all $x,x'\in X$.
Denote by $$\Lip(f)=\sup_{x\ne 
x'}\{\frac{d_Y(f(x),f(x'))}{d_X(x,x')}\}$$
the minimal Lipschitz constant of $f$.

\mk\noindent Every simplicial complex $K$ carries a metric
where
all simplexes are isometric to the standard euclidean simplex.
We will call the maximal such metric on $K$ {\it uniform} and usually we will
denote the corresponding metric space by $K_U$. \footnote{An important 
exception to this convention occurs with {\it asymptotic polyhedra}, 
which are defined below.} Note that the metric space $K_U$ is geodesic.
If no metric is specified,
we will assume that a finite complex is supplied with this
uniform metric. In particular, the complexes in the two lemmas below 
are assumed to have the uniform metric.

\mk\noindent The following lemma is a special case of Theorem A from [SW].
\proclaim{Lemma 3.1} Let $Y$ be a finite simplicial complex with
$\pi_n(Y)$ finite. Then for every $\lambda>0$ there is a $\mu>0$ such
that every map $f:B^n\to Y$ with $\Lip(f|_{S^{n-1}})\le\lambda$ can
be deformed to a $\mu$-Lipschitz map $g:B^n\to Y$ by means of a
homotopy $h_t:B^n\to Y$ with $h_t|_{S^{n-1}}=f|_{S^{n-1}}$.
\endproclaim

\proclaim{Lemma 3.2} Let $L$ be a finite dimensional complex and
let $K$ be a finite complex with finite homotopy groups $\pi_i(K)$
for $i\le\dim L+1$. Let $f,g:L\to K$ be homotopic Lipschitz maps.
Then every homotopy between $f$ and $g$ can be deformed to a
Lipschitz homotopy $H:L\times[0,1]\to K$.
\endproclaim
\demo{Proof} Let $F:L\times I\to K$ be a homotopy between $f$ and
$g$. By induction on $n$ using Lemma 3.1, we construct a
$\mu_n$-Lipschitz map $H_n:L^{(n)}\times I\cup L\times\{0,1\}\to
K$ which is a deformation of $F$ restricted to the $n$-skeleton
$L^{(n)}$ such that $H_n$ extends $H_{n-1}$. Here the fact that
$L$ (and hence $L\times I$) is geodesic is essential for the
argument because this condition guarantees that the union of
$\lambda$-Lipschitz maps on $\Delta^n\times I$ is
$\lambda$-Lipschitz on $L^{(n)}\times I$ (see [Dr2] or [Dr3] for
details). Then $H=H_m$ for $m=\dim L$. \qed
\enddemo

\mk\noindent Let $x_0\in X$ be a basepoint in a metric space $X$. For $x\in X$ and
$A\subset X$ we denote by $\|x\|=\dist(x,x_0)$ and $\|A\|=\max\{\|z\|\,|\,z\in A\}$.
By $B_r(x)$ we denote the closed $r$-ball in $X$ centered at $x$. We use notation
$B_r=B_r(x_0)$.

\mk\noindent Let $\sigma$ be an $n$-dimensional simplex spanned in a Euclidean space.
By $s_{\sigma}:\sigma\to\Delta^n$
we denote a simplicial homeomorphism onto the standard $n$-simplex
$\Delta^n=\{(x_i)\in\R^{n+1}\mid \sum x_i=1, x_i\ge0\}$.
A locally finite simplicial complex $L$ with a geodesic metric on it is called
an {\it asymptotic polyhedron} if
every simplex in $L$ is isometric to a simplex $\sigma$ spanned in a
Euclidean space and
$\lim_{\|\sigma\|\to\infty}\Lip(s_{\sigma})=0$. Asymptotic polyhedra 
usually do not have the uniform metric.

\proclaim{Proposition 3.3} Let $(K,d)$ be an asymptotic polyhedron
and let $d_U$ denote the uniform geodesic metric on $K$. Then the
identity map $u:(K,d)\to (K,d_U)$ satisfies the condition $\lim_{x 
\to \infty} \diam 
(u(B_{r}(x))) = 0$.
\endproclaim
\demo{Proof} Let $r$ and $\epsilon>0$ be given. There is a $t$ so 
that $\Lip(s_{\sigma})<\epsilon/(2r)$ for every $\sigma$ with $\|\sigma\|\ge
t$. Choose $x$ so that $\|x\|\ge t+2r$. For every two points $z,z'\in
B_r(x)$, a geodesic segment $J=[z,z']$ joining them does not
intersect $B_t$. There is a partition of $J$:
$z=z_0<z_1\dots<z_m=z'$ such that every segment
$J_i=[z_i,z_{i+1}]$ lies in a simplex $\sigma_i$. By our choice of $t$,
$\Lip(s_{\sigma_i})<\epsilon/(2r)$ for each $i$. Therefore
$d_U(z_i,z_{i+1})\le\epsilon/(2r)d(z_i,z_{i+1})$ and hence
$d_U(z,z')\le\sum\epsilon/(2r)d(z_i,z_{i+1})=\epsilon/(2r)d(z,z')\le\epsilon$.\qed
\enddemo

\mk\noindent A family of subsets $\sA$ is said to be $d$-disjoint if for 
every pair of 
sets $A$ and $A'$ in $\sA$ we have $d(x,x')>d$ for all $x\in A$ and 
$x'\in A'$. A family of subsets $\sU$ of a metric space $X$ is {\it 
uniformly bounded} if there is a $B>0$ so that the diameter of each 
element of $\sU$ is less than $B$. 
 Gromov,  [G1], defined the {\it asymptotic dimension} $\as X$ of a metric space $X$ as follows:
$\as X\le n$ if for every $d$ there are $n+1$ $d$-disjoint uniformly bounded families
$\sU_i$, $i=0,\dots, n$ of subsets of $X$ such that the union $\sU=\cup\sU_i$ is a cover
of $X$. 

\mk\noindent Gromov's definition can be equivalently reformulated as follows: $\as X\le n$
if for any arbitrary large number $d$
there is a uniformly bounded open cover $\sU$ of $X$ with multiplicity $\le n+1$
and with Lebesgue number $\ge d$ (see Assertion 1 in [BD2] for a proof).

\mk\noindent We note that for every $n$-dimensional asymptotic polyhedron $L$,
$\as L \le n$.

\mk\noindent Let $\sU$ be a locally finite open cover of a metric space $X$. The {\it canonical
projection} to the nerve $p:X\to \Nerve(\sU)$ is defined using the   
partition of unity $\{\phi_U:X\to\R\}_{U\in\sU}$ defined by 
$\phi_U(x)=d(x,X\setminus U)/\sum_{V\in\sU} d(x,X\setminus V)$.
The family $\{\phi_U:X\to\R\}_{U\in\sU}$ defines a map $p$ to the
Hilbert space $l_2(\sU)$ of  square summable functions on $\sU$ with the Dirac functions
$\delta_U$, $U\in \sU$ as the basis. The nerve
$N(\sU)$ of the cover $\sU$ is realized in $l_2(\sU)$ by
taking every vertex $U$ to $\delta_U$.
Clearly, the image of $p$ lies in the nerve.
Given a family of positive numbers $\bar\lambda=\{\lambda_U\}_{U\in\sU}$ we can change the above
imbedding of the nerve $N(\sU)$ into $l_2(\sU)$ by taking each vertex $U$ to $\lambda_U\delta_U$.
Then the projection $p^{\bar\lambda}_{\sU}:X\to N(\sU)$ to modified realization is given by the formula
$p^{\bar\lambda}_{\sU}(x)=\sum_{U\in\sU}\lambda_U\phi_U(x)\delta_U$.

\mk\noindent For a subset $A\subset X$ we denote by $L(\sU)|_A$ the Lebesgue
number of $\sU$ restricted to $A$. More precisely, $L(\sU)|_A=\inf_{y\in
A}\max_{V\in\sU}d(y,X\setminus V)$.

\proclaim{Proposition 3.4} Let $\sU$ be a locally finite cover of
a geodesic metric space $X$ with the multiplicity $\le m$. Then
the above projection to the nerve $p^{\bar\lambda}_{\sU}:X\to
N=N(\sU)$ for $\bar\lambda=\{\lambda_U\}$ with
$\lambda_U=L(\sU)|_U/(2m+1)^2$ is 1-Lipschitz
 where the nerve is taken with the
intrinsic metric induced from $l_2(\sU)$.
\endproclaim
\demo{Proof} We will show that the map $\bar p=p^{\bar\lambda}_{\sU}$
is 1-Lipschitz as a map to $l_2(\sU)$. Then the partition of $N$ into
simplices defines a locally finite partition on $X$ such that
$\bar p$ is 1-Lipschitz on every piece, considered as the map to
$N$ with the induced path metric. Since $X$ is a geodesic metric
space, this will imply that $\bar p$ is 1-Lipschitz.

\mk\noindent Let $x,y\in X$ and $U\in\sU$.
The triangle inequality implies
$$|d(x,X\setminus U)-d(y,X\setminus U)|\le d(x,y).$$
It is easy to see that $\sum_{V\in\sU}d(x,X\setminus V)\ge L(\sU)|_U=(2m+1)^2\lambda_U$ for $x\in U$. Then
$|\phi_U(x)-\phi_U(y)|= |\frac{d(x,\,X \setminus U)}{\sum_{V \in 
\sU}d(x,\,X\setminus V)}-\frac{d(y,\,X \setminus U)}{\sum_{V \in 
\sU}d(x,\,X\setminus V)}+\frac{d(y,\,X \setminus U)}{\sum_{V \in 
\sU}d(x,\,X\setminus V)}-\frac{d(y,\,X \setminus U)}{\sum_{V \in 
\sU}d(y,\,X\setminus V)}|\le$
$$\frac{1}{\sum_{V\in\sU}d(x,X\setminus V)}d(x,y)+
d(y,X\setminus U)|\frac{1}{\sum_{V\in\sU}d(x,X\setminus V)}-
\frac{1}{\sum_{V\in\sU}d(y,X\setminus V)}|$$
$$\le\frac{1}{(2m+1)^2\lambda_U}d(x,y)+
\frac{1}{(2m+1)^2\lambda_U}(\sum_{V\in\sU}|d(x,X\setminus V)-d(y,X\setminus V)|)\le
\frac{1}{\sqrt{2m+1}\lambda_U}d(x,y).$$
Then
$\|p(x)-p(y)\|=$
$$=(\sum_{U\in\sU}\lambda_U^2(\phi_U(x)-\phi_U(y))^2)^{\frac{1}{2}}\le
((2m)\frac{1}{2m+1}d(x,y)^2)^{\frac{1}{2}}
\le d(x,y).$$\qed
\enddemo

\proclaim{Lemma 3.5} Suppose that $X$ is a geodesic metric space
with $\as X\le n$ and let $f:X\to\R_+$ be a proper function.
Then there exist a compact set $C\subset X$ and a map $\phi:X\to N$
to a $2n+1$-dimensional asymptotic polyhedron with
$\diam(\phi^{-1}(\Delta))\le f(z)$ for all
$z\in\phi^{-1}(\Delta)\setminus C$ and for all simplices $\Delta$ of 
$N$. Moreover, every vertex of $N$ lies in the image of $\phi$.
\endproclaim
 \demo{Proof}
We construct $\phi$ as the projection $p_{\sU}^{\bar\lambda}$
from Proposition 3.4
to the nerve of a cover $\sU$ of $X$.

\mk\noindent Fix a monotone sequence $\{l_i\}$ tending to infinity. Since $\as X\le n$,
for every
$i$ there is a uniformly bounded cover $\sU_i$ of multiplicity
$n+1$ with the Lebesgue number $L(\sU_i)>l_i$. Let $m_i$ be an
upper bound for the diameter of elements of $\sU_i$. Choose a
sequence $\{r_i\}$ such that  $f(X\setminus B_{r_{i}-m_{i}}) \ge
2m_{i+1}$ and take $\sU=\cup_i\sU'_i$ where
$$\sU'_i=\{U\setminus B_{r_{i-1}}\mid U\in\sU_i, U\cap B_{r_i}\ne\emptyset\}.$$

\mk\noindent We take $C=B_{r_1}$ and check that the conclusion of our Lemma holds. By 
definition, the preimage $\phi^{-1}(\Delta)$ lies in a union $
\cup_{j \in J} U'_j$ of sets from $\sU$ with nonempty intersection. Let
$x\in\cap_{j \in J} U'_j$. Each $U'_j$ belongs to the cover $\sU'_i$ for some
$i$. Let $k$ be maximal among those $i$'s for all $j \in J$. Then $x\notin
B_{r_{k-1}}$. Therefore $U'_j\subset X\setminus
B_{r_{k-1}-m_{k-1}}$ for all $j$. Hence $f(z)\ge 2m_k\ge
\diam\cup_jU'_j$ for $z\in X\setminus B_{r_{k-1}-m_{k-1}}$ and
$k>1$.

\mk\noindent Since $l_i\to\infty$ the nerve $N$ realized in $l_2(\sU)$ as above with the
intrinsic metric is an asymptotic
polyhedron. Since the covers $\sU'_i$ and $\sU'_j$ do not intersect for $|i-j|>1$,
the multiplicity of $\sU$ is $\le 2n+2$.

\mk\noindent Finally, if we replace $\sU$ by an irreducible 
subcover, $\phi$ is onto the vertices of $N$.
\qed
\enddemo
\mk\noindent REMARK. By a standard dimension theoretic trick the polyhedron $N$
can be chosen to be $n$-dimensional.

\mk\noindent A metric space $X$ is called {\it uniformly contractible} if there is
a function $S:\R_+\to\R_+$ such that every ball $B_r(x)$ is contractible to
a point in the ball $B_{S(r)}(x)$. The function $S$ is called a {\it 
contractibility function} for $X$. We will always take our 
contractibility functions to be strictly monotone.

\proclaim{Lemma 3.6} Let $X$ be a uniformly contractible proper
metric space with $\as X\le n$. Then given a proper function
$g:X\to\R_+$ there exist a $2n+1$-dimensional asymptotic polyhedron
$N$, a proper 1-Lipschitz map $\phi:X\to N$, and a proper homotopy
inverse map $\gamma:N\to X$ with $d(x,\gamma\phi(x))<g(x)$ for all
$x\in X$. Moreover, there is a compact set $C\subset X$ such that
$\diam(\phi^{-1}(\Delta))\le g(x)$ for all $x\in\Delta\setminus C$.
\endproclaim
\demo{Proof} Let $S$ be a contractibility function. We define
$\rho(t)=S^{-1}(t/2)$ where $S^{-1}$ is the inverse function for
$S$. Then we take $f=\rho^{2n+1}\circ g$, the composition of $ g$
and the $2n+1$-fold iteration of $\rho$. Clearly $g\le f$. We assume
here that $g(x)\le\|x\|/2$. We define by induction on $i$ a lift $\gamma$ on the
$i$-skeleton $N^{(i)}$ of the nerve of a cover of $X$ given by
Lemma 3.5 for this choice of $f$. We take
$\gamma(v)\in\phi^{-1}(v)$ for every vertex $v$.  Then using the
uniform contractibility of $X$ we can extend $\gamma$ with control
over the 1-skeleton $N^{(1)}$ and so on.  Without loss of generality
(and for simplicity of exposition) we may assume that $X$ is a
polyhedron of the dimension $n$ supplied with a triangulation of mesh
$\le 1$.  By induction on $i$ we define a homotopy $H:X^{(i)}\times
I\to X$ joining the identity map with $\gamma\circ\phi$.  We consider
a function $\psi(x)=\|x\|-\max\{d(x,y)\mid y\in H(x\times I)\}$.  If
$\psi$ tends to infinity, then the map $H$ is proper, but it is easy to verify that $\psi(x)\ge \|x\|/2$.
 \qed
\enddemo

\mk\noindent We recall that the Higson compactification
$\bar X$ of a proper metric space $X$ can be defined as the
maximal ideal space of the completion of the ring of bounded
functions with the gradient tending to zero at infinity [Ro1]. The
defining property of the Higson corona is the following:

\mk\noindent(*){\it  A continuous map $f:X\to Z$ to a compact
metric space is extendable to the Higson corona $\nu X$
if and only if it satisfies the condition: For arbitrary large} $R$
$$\lim_{\|x\|\to\infty}\diam(f(B_R(x))=0.$$

\mk\noindent Note that a proper Lipschitz map $f:X\to Y$ induces
a continuous mapping between the Higson coronas $\bar f:\nu X\to\nu Y$.

\proclaim{Theorem 3.7} Let $X$ be a uniformly contractible
geodesic proper metric space with finite asymptotic dimension
and let $\bar X$ be the Higson compactification. Then
 $\check H^{n}(\bar X;\Z_p)=0$ for all $n>0$ and all $p$.
\endproclaim
\demo{Proof}
We show that every map $\alpha:\bar X\to K(\Z_p,n)$ is null homotopic.
Since $\bar X$ is compact, the image $\alpha(\bar X)$ is contained in the
$k$-skeleton $K=K(\Z_p,n)^{(k)}$ for some $k$.
We will assume that $K$ is a finite complex of dimension at least $2n+2$. 
For convenience, we replace $K$ by a Riemannian manifold with 
Riemannianly collared boundary. Let $4\epsilon_K$ be a convexity radius 
for $K$, i.e. a positive constant such that every two points
in an $4\epsilon_K$-ball can be joined by a unique geodesic in that 
$4\epsilon_{K}$-ball.
 Since the map $\alpha|_X\to K$ is extendable over the Higson corona,
the function $R_{\alpha}(t)=\Lip(\alpha|_{X\setminus B_t(x_0)})$
tends to zero at infinity. We apply Lemma 3.6 with
$g(x)\le\min\{\epsilon_K/R_{\alpha}(\|x\|/2), \|x\|/4\}$ to obtain
an asymptotic polyhedron $N$ and maps $\phi:X\to N$ and
$\gamma:N\to X$. We may assume that $\phi$ is surjective on vertices. Let
$[u,v]$ be an edge in $N$, then
$d_K(\alpha\gamma(u),\alpha\gamma(v)) \le
R_{\alpha}(t_0)d_X(\gamma(u),\gamma(v))$ where
$t_0=\min\{\|\gamma(u)\|,\|\gamma(v)\|\}$. We may assume that
there are $x,y\in X$ such that $\phi(x)=u$ and $\phi(y)=v$. Then
$$d_X(\gamma(u),\gamma(v))\le d_X(x,y)+
\epsilon_K/R_{\alpha}(\frac{1}{2}\|x\|)+\epsilon_K/R_{\alpha}(\frac{1}{2}\|y\|)\le
\diam\phi^{-1}[u,v]+2\epsilon_K/R_{\alpha}(\frac{1}{2}\|x\|)$$
provided that $\|x\|\le\|y\|$. Because of the inequality
$g(x)\le\|x\|/4$ we have that $R_{\alpha}(t_0)\le
R_{\alpha}(\|x\|/2)$. By Lemma 3.6, 
$\diam(\phi^{-1}([u,v]))<g(x)$. As a result, we obtain the
inequality $d_K(\alpha\gamma(u),\alpha\gamma(v))\le 3\epsilon_K$.
This means that there is a map $\beta:N \to K$ which coincides on vertices 
with the map $\alpha\circ\gamma$ and which is
$3\epsilon_K$-Lipschitz for $N_{U}=N$ taken with the uniform metric: 
$\beta$ is obtained by linear extension of the restriction of  
$\alpha\circ\gamma$ to vertices. Note that $\beta$ is 
$\epsilon_{K}$-close to $\alpha\circ\gamma$ and that these maps are 
therefore homotopic.

\mk\noindent Since $X$ is contractible, the map $\alpha\circ\gamma$ is null
homotopic. Therefore, so is $\beta$. Note that the homotopy groups $\pi_i(K)$ are finite for
$i\le\dim N+1$. We apply Lemma 3.2 to obtain a $\lambda$-Lipschitz
homotopy $H:N_U\times I\to K$ of $\beta$ to a constant
map. This homotopy defines a Lipschitz map $\tilde H:N_U\to
K^I_{\lambda}$ to the space of $\lambda$-Lipschitz mappings of the
unit interval $I$ to $K$.
We note that the space $K^I_{\lambda}$ is compact. Then, by Proposition 3.3,
$\tilde H\circ u:N\to K^I_{\lambda}$ satisfies the Higson extendibility
condition (*).
Let $\bar h:\bar N\to K^I_{\lambda}$ be the extension over the Higson corona.
This extension defines a map $\bar H:\bar N\times I\to K$. The map $\bar H$ is
a homotopy between the extension $\overline{\beta}$ and a constant map.
To complete the proof, we show that $\alpha$ is homotopic to $\bar\beta\circ\bar\phi$
where $\bar\phi$ is the extension of the Lipschitz map $\phi$ to the Higson
compactifications. Note that
$$d_K(\alpha(x),\alpha\gamma\phi(x))\le R_{\alpha}(t_0)d(x,\gamma\phi(x))\le
R_{\alpha}(t_0)\epsilon_K/R_{\alpha}(\|x\|/2)\le\epsilon_K$$ where
$t_0=\min\{\|x\|,\|\gamma\phi(x)\|\}\ge\|x\|/2$. Therefore, 
$d_{K}(\alpha(x),\,\beta\circ \phi(x))\le 4\epsilon_{K}$. Then for every $x\in X$ we join
the points $\alpha(x)$ and $\beta\phi(x)$ by the unique geodesic
$\psi_x:I\to K$.
This defines a map $\tilde\psi:X\to K^I_{\mu}$. Since both
$\alpha$ and $\beta\circ\phi$ satisfy the condition (*), the map
$\tilde\psi$ has the property (*). Let $\bar\psi:\bar X\to K^I_{\mu}$ be the
extension of $\tilde\psi$ to the Higson corona. The map $\bar\psi$ defines
a homotopy $\Psi:\bar X\times I\to K$ between $\alpha$ and $\bar\beta\circ\bar\phi$.
\qed
\enddemo

\mk\noindent We have now proven Theorem 3. The following asserts 
Theorem 2:

\proclaim{Corollary 3.8} Suppose that a group $\Gamma$ has 
finite asymptotic dimension and that $\BG$ is a finite complex.
Then the integral assembly map $A:H_*(\BG;\L)\to
L_*(\Z\Gamma)$ is a monomorphism.
\endproclaim
\demo{Proof}
In view of the inequality $\dim\nu\Gamma\le\as \Gamma$ [DKU] we can apply
Theorem 2.4. \qed
\enddemo
\proclaim{Corollary 3.9} If $\Gamma$ and $R$ satisfy the finiteness 
conditions of Theorem 2, then $H_{i}(B\Gamma;\,\L(R)) \to 
L_{i}(R\Gamma)$ is injective for all $i$.
\endproclaim

\enddocument